\documentclass{article}

\usepackage{amsmath}
\usepackage{amssymb}
\usepackage{amsthm}
\usepackage{amsfonts}
\usepackage{latexsym}
\usepackage[pdftex]{hyperref}
\usepackage{epsfig}

\usepackage{MnSymbol}

\newtheorem{Theorem}{Theorem}[section]
\newtheorem{Lemma}[Theorem]{Lemma}
\newtheorem{Corollary}[Theorem]{Corollary}
\newtheorem{Definition}[Theorem]{Definition}

\newtheorem*{Theorem*}{Theorem}

\newcommand{\G}{\mathbb{G}}

\title{Seifert-Tait graphs}
\author{Stephen Huggett and Alina Vdovina}

\begin{document}
\maketitle
\begin{abstract}
We show that among alternating knots, those which have diagrams whose Seifert and Tait graphs are isomorphic are dominant.
\end{abstract}

\section{Introduction}

A knot diagram $D$ defines two graphs, the \emph{Seifert} graph and the \emph{Tait} graph. The Seifert graph arises from the (canonical) Seifert surface of $D$, while the Tait graph arises from a chessboard colouring of $D$. Evidently, these procedures are quite different, and one might expect them to lead to different graphs, which they often do. However, it is not difficult to find knot diagrams whose Seifert and Tait graphs are isomorphic.

\begin{Definition}
A \emph{Seifert-Tait knot} is a knot having an alternating diagram whose Seifert and Tait graphs are isomorphic.
\end{Definition}

Our main theorem is to show that Seifert-Tait knots dominate the alternating knots:

\begin{Theorem*}
For a fixed genus $g>1$,
$$\frac{\#\{\hbox{Seifert-Tait alternating knots with crossing number }n\hbox{ and genus }g\}}{\#\{\hbox{prime alternating knots with crossing number }n\hbox{ and genus }g\}}\longrightarrow 1$$
as $n\longrightarrow\infty$.
\end{Theorem*}

Here is an outline of the paper. Section 2 is a review of the definitions of Seifert graphs and Tait graphs, and the relationship between them demonstrated in \cite{HMV}. In Section 3 we define \emph{flat} diagrams, and show that they must have isomorphic Seifert and Tait graphs. Section 4 focusses on the \emph{the cyclic word algorithm} introduced in \cite{Alina2005}, different from Seifert's algorithm, for constructing the Seifert surface of a knot. This completes the groundwork for Section 5, in which we review the argument given in \cite{Alina2005} proving that among alternating knots, those having flat diagrams dominate. Then we deduce our main theorem.

\section{Seifert and Tait graphs}

There is an algorithm due to Seifert which, given an oriented knot diagram $D$, constructs an orientable surface---presented as a ribbon graph---whose boundary is that knot. We call this surface the \emph{canonical Seifert surface} of the knot. The algorithm is well known, of course, but we will be referring to it in detail, so we review it here to set the terminology.

\begin{figure}[h]
\centering
\includegraphics[trim= 2cm 1.3cm 2cm 1cm, clip=true]{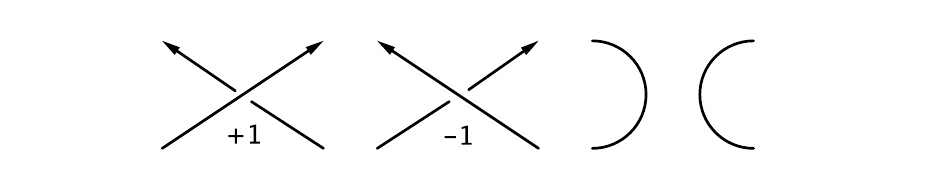}
\caption{Signs of oriented crossings, and the only possible splicing.}
\label{Writhe}
\end{figure}

Replace each crossing as indicated in figure \ref{Writhe}: this is known as \emph{splicing}. The diagram becomes a set of disjoint simple closed curves in the plane. These closed curves are called \emph{Seifert circles}. They may be nested: let $h(\gamma)$ be the number of Seifert circles that contain the Seifert circle $\gamma$. For each Seifert circle $\gamma$, take a disc at height $h(\gamma)$ above the plane whose boundary projects down onto $\gamma$. These discs inherit an orientation from the diagram. To construct the required canonical Seifert surface, insert a half-twisted ribbon at the site of each crossing, the sign of the twist being determined by the sense of the crossing. The resulting surface is exhibited as a ribbon graph having the original knot as its only boundary component. See \cite{Cromwell2004} for a more detailed, and beautifully clear, discussion.

From the canonical Seifert surface we may define an (abstract) graph as follows. Each disc of the canonical Seifert surface is a vertex of the graph, and each ribbon of the canonical Seifert surface is an edge of the graph. This graph is called the \emph{Seifert graph} of $D$.

\begin{Lemma}
Seifert graphs are planar and bipartite.
\end{Lemma}

Note that we shall be regarding a \emph{planar} graph as having an embedding in the \emph{sphere} rather than the plane.

\begin{figure}[h]
\centering
\includegraphics[height=2cm]{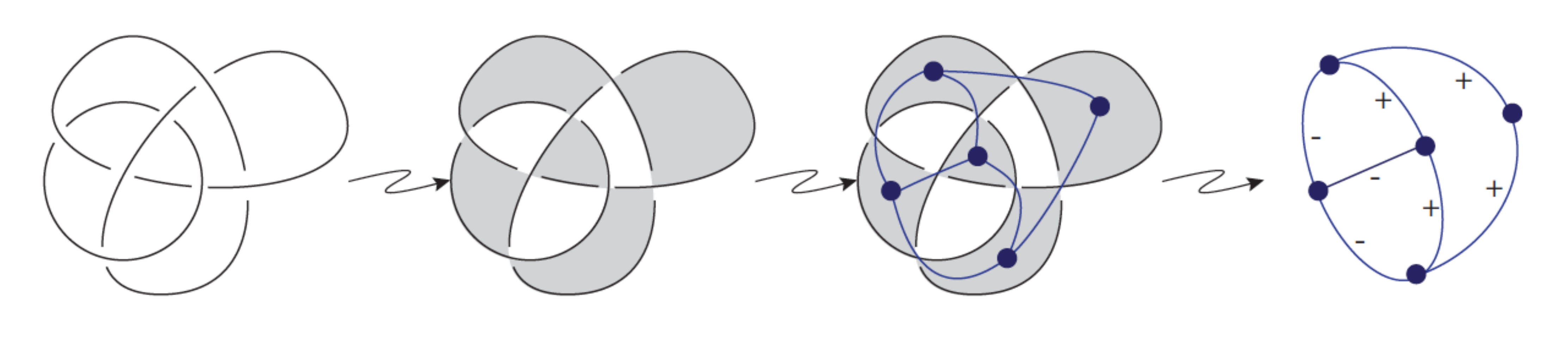}
\caption{Moving between a link diagram and its signed Tait graph.}
\label{SignedTaitExample}
\end{figure}

Returning to our knot diagram $D$, we may define another planar graph, called its \emph{Tait graph} $T$. We colour the components of the complement of $D$ black or white to form a ``chessboard" pattern. The black components are the vertices of $T$, and the crossings of $D$ are the edges of $T$. See figure \ref{SignedTaitExample}. We could have chosen the white components as the vertices instead, when we would have obtained the plane dual graph $T^{*}$. In other words, the knot diagram $D$ actually gives rise to \emph{two} Tait graphs, $T$ and $T^{*}$.

The Seifert graph may be described in terms of the two Tait graphs as follows. See \cite{HMV} and \cite{Bipartite}.

\setlength{\unitlength}{1pt}
\begin{figure}
\begin{picture}(150,110)(-100,-20)

\put(13,55){\circle*{4}} \put(77,55){\circle*{4}}
\put(13,55){\line(1,0){64}} \qbezier(10,65)(35,65)(45,55)
\qbezier(45,55)(55,45)(80,45) \qbezier(10,45)(35,45)(43,53)
\qbezier(47,57)(55,65)(80,65) \put(15,65){\vector(1,0){0}}
\put(11,45){\vector(-1,0){0}} \put(76,65){\vector(-1,0){0}}
\put(80,45){\vector(1,0){0}}

\put(133,55){\circle*{4}} \put(197,55){\circle*{4}}
\put(133,55){\line(1,0){64}} \qbezier(130,65)(155,65)(165,55)
\qbezier(165,55)(175,45)(200,45) \qbezier(130,45)(155,45)(163,53)
\qbezier(167,57)(175,65)(200,65) \put(135,65){\vector(1,0){0}}
\put(135,45){\vector(1,0){0}} \put(200,65){\vector(1,0){0}}
\put(200,45){\vector(1,0){0}}

\put(13,5){\circle*{4}} \put(77,5){\circle*{4}}
\put(13,5){\line(1,0){64}} \qbezier(10,15)(40,15)(40,5)
\qbezier(50,5)(50,-5)(80,-5) \qbezier(10,-5)(40,-5)(40,5)
\qbezier(50,5)(50,15)(80,15) \put(43,-10){d}
\put(15,15){\vector(1,0){0}} \put(11,-5){\vector(-1,0){0}}
\put(76,15){\vector(-1,0){0}} \put(80,-5){\vector(1,0){0}}

\put(133,5){\circle*{4}} \put(197,5){\circle*{4}}
\put(133,5){\line(1,0){64}} \qbezier(130,15)(165,5)(200,15)
\qbezier(130,-5)(165,5)(200,-5) \put(165,-10){c}
\put(134,14){\vector(3,-1){0}} \put(134,-4){\vector(3,1){0}}
\put(200,15){\vector(3,1){0}} \put(200,-5){\vector(3,-1){0}}

\put(25,75){d-edge} \put(145,75){c-edge}
\put(-60,50){$\textrm{before}$} \put(-60,0){$\textrm{after}$}

\end{picture}
\caption{d-edges and c-edges before and after splicing a
crossing.}
\label{c and d edges}
\end{figure}
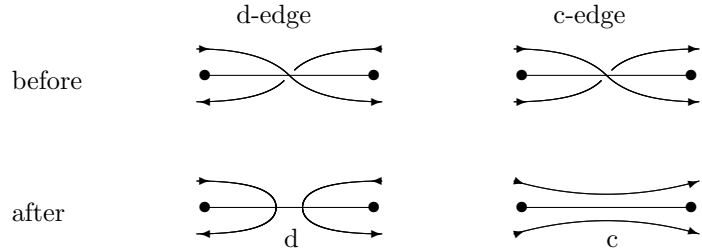

In the Seifert algorithm, each edge of $T$ is either contracted (and then called a \emph{c-edge}) or deleted (and then called a \emph{d-edge}), as shown in figure \ref{c and d edges}. Let $E_{c}$ be the set of c-edges of $T$, and let $V_{c}$ be the set of vertices of $T$ adjacent to a c-edge. Let $C=(V_{c},E_{c})$. The graph $C$ is not necessarily connected, but each component is Eulerian. The d-edges in $T$ correspond to c-edges in $T^{*}$. Let $E'_{c}$ be the set of c-edges of $T^{*}$, and let $V'_{c}$ be the set of vertices of $T^{*}$ adjacent to a c-edge. Let $C'=(V'_{c},E'_{c})$. The graph $C'$ is not necessarily connected, but each component is Eulerian.

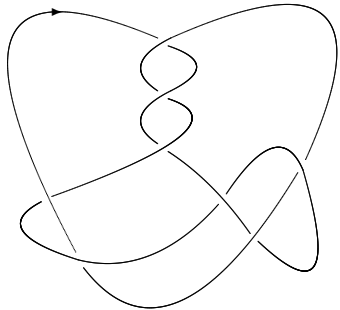
\begin{figure}[h]
\begin{picture}(300,120)(-130,-20)

\qbezier(14,-1)(-50,120)(45,80)
\put(9,90){\vector(1,0){0}}
\qbezier(49,77)(70,70)(50,60)
\qbezier(50,60)(30,50)(45,40)
\qbezier(49,37)(66,26)(80,6)
\qbezier(100,30)(117,-29)(83,3)
\qbezier(71,21)(91,51)(100,30)
\qbezier(15,-4)(43,-11)(68,17)
\qbezier(1,18)(-20,7)(15,-4)
\qbezier(49,57)(80,47)(5,20)
\qbezier(50,80)(30,70)(45,60)
\qbezier(50,80)(140,120)(101,34)
\qbezier(17,-7)(50,-50)(98,29)

\end{picture}
\caption{A diagram of the knot $8_{11}$, with an orientation.}
\label{8_{11}}
\end{figure}

\begin{figure}[h]
\begin{picture}(150,70)(-70,40)

\put(-150,30){\begin{picture}(100,100)(-70,40)
\put(50,110){\circle*{5}}
\put(150,110){\circle*{5}}
\qbezier{(50,110),(100,125),(150,110)}
\put(50,110){\line(50,0){100}}
\qbezier{(50,110),(100,95),(150,110)}

\put(50,60){\circle*{5}}
\put(150,60){\circle*{5}}
\put(50,60){\line(50,0){100}}
\put(50,110){\line(0,-50){50}}
\put(150,110){\line(0,-50){50}}

\put(175,85){\circle*{5}}
\put(150,110){\line(1,-1){25}}
\put(150,60){\line(1,1){25}}

\put(100,112){c}
\put(100,104){c}
\put(100,97){c}

\put(100,55){c}
\put(51,83){c}
\put(151,83){c}

\end{picture}}

\put(0,30){\begin{picture}(100,100)(-90,40)

\put(50,110){\circle*{5}}
\put(50,90){\circle*{5}}
\put(50,60){\circle*{5}}
\put(50,110){\line(0,-50){20}}
\put(50,90){\line(0,-50){30}}

\put(100,110){\circle*{5}}
\put(100,60){\circle*{5}}
\put(50,110){\line(50,0){50}}
\put(50,60){\line(50,0){50}}

\qbezier{(50,60),(71,90),(100,110)}
\qbezier{(50,60),(79,80),(100,110)}

\qbezier{(100,60),(93,80),(100,110)}
\qbezier{(100,60),(107,80),(100,110)}

\put(97,81){c}
\put(104,81){c}

\end{picture}}
\end{picture}
\caption{The Tait graphs $T$ and $T^{*}$. The unlabelled edges are all d-edges.}
\end{figure}
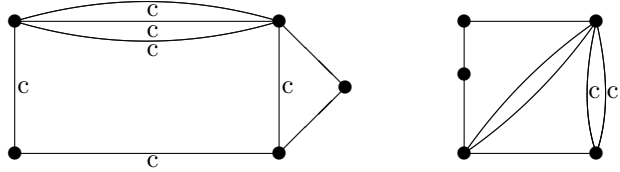

\begin{figure}[h]
\begin{picture}(150,80)(-70,40)

\put(-90,30){\begin{picture}(100,100)(-70,40)
\put(50,110){\circle*{5}}
\put(150,110){\circle*{5}}
\qbezier{(50,110),(100,125),(150,110)}
\put(50,110){\line(50,0){100}}
\qbezier{(50,110),(100,95),(150,110)}

\put(50,60){\circle*{5}}
\put(150,60){\circle*{5}}
\put(50,60){\line(50,0){100}}
\put(50,110){\line(0,-50){50}}
\put(150,110){\line(0,-50){50}}

\end{picture}}

\put(0,30){\begin{picture}(100,100)(-70,40)

\put(100,110){\circle*{5}}
\put(100,60){\circle*{5}}

\qbezier{(100,60),(93,80),(100,110)}
\qbezier{(100,60),(107,80),(100,110)}

\end{picture}}
\end{picture}
\caption{The graph $\Phi = i(C) \cup i(C')$.}
\label{phigraph}
\end{figure}
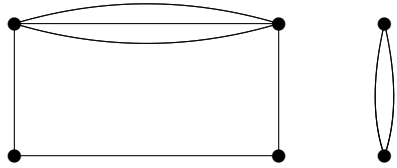

\begin{figure}[h]
\begin{picture}(150,50)(-70,40)

\put(50,60){\circle*{5}}
\put(150,60){\circle*{5}}
\qbezier{(50,60),(100,75),(150,60)}
\put(50,60){\line(50,0){100}}
\qbezier{(50,60),(100,45),(150,60)}

\put(50,110){\circle*{5}}
\put(150,110){\circle*{5}}
\put(50,110){\line(50,0){100}}
\put(50,110){\line(0,-50){50}}
\put(150,110){\line(0,-50){50}}

\put(250,60){\circle*{5}}
\qbezier{(150,60),(200,75),(250,60)}
\qbezier{(150,60),(200,45),(250,60)}

\end{picture}
\caption{The Seifert graph $S=\Phi^{*}$.}
\end{figure}
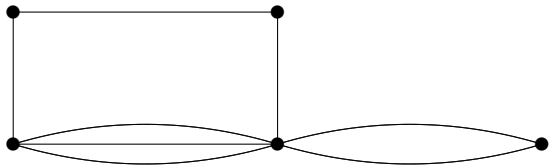

The Tait graph $T$ is equipped with an embedding $i:T\rightarrow S^{2}$. This induces an embedding of its dual $T^{*}\rightarrow S^{2}$, which we also denote by $i$. Now take $i(C)\cup i(C')$, and denote the resulting graph by $\Phi$. See figure \ref{phigraph}. $\Phi$ is embedded in $S^{2}$, but not necessarily cellularly.

\begin{Theorem}[\cite{HMV}]
$\Phi^{*}=[i(C)\cup i(C')]^{*}$ is the Seifert graph of $D$.
\end{Theorem}

\section{Flat diagrams}

Recall that we are regarding a planar graph as having an embedding in the sphere.

\begin{Definition}
A Seifert circle is of \emph{type I} if it bounds a disc which does not contain a crossing of the knot diagram. Otherwise it is of \emph{type II}.
\end{Definition}

\begin{Definition}[\cite{Alina2005}]
A \emph{flat} diagram is a knot diagram all of whose Seifert circles are of type I.
\end{Definition}

In the literature (\cite{Alina2005}, \cite{Cromwell1989}, and \cite{Murasugi}, for example) these are referred to as \emph{special} diagrams. But it is shown in \cite{Alina2005} that for a fixed genus $g>1$ the special alternating knots dominate the prime alternating knots (see Theorem \ref{limit}), so in that sense they are not ``special". We have chosen the word \emph{flat}, because if all the Seifert circles are of type I then a projection from $S^{2}$ to the plane can be chosen so that all the Seifert circles have height $0$. In this case the Seifert graph is not just planar, it is \emph{already} embedded in $S^{2}$.

\begin{Lemma}[\cite{Cromwell2004}]
\label{notblock}
If a knot diagram is not flat then its Seifert graph has a cut vertex.
\end{Lemma}
As stated, the converse is false: see Figure \ref{FlatCutVertex} for a counterexample. However, this comes from a knot diagram with a nugatory crossing, which may be reduced to the diagram of a non-prime knot.

\setlength{\unitlength}{0.5pt}
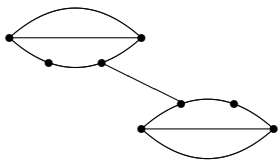
\begin{figure}[h]
\begin{picture}(150,150)(-150,0)

\put(50,110){\circle*{5}}
\put(150,110){\circle*{5}}
\qbezier{(50,110),(100,155),(150,110)}
\put(50,110){\line(50,0){100}}
\qbezier{(50,110),(100,65),(150,110)}
\put(80,91){\circle*{5}}
\put(120,91){\circle*{5}}

\put(150,41){\circle*{5}}
\put(250,41){\circle*{5}}
\qbezier{(150,41),(200,86),(250,41)}
\put(150,41){\line(50,0){100}}
\qbezier{(150,41),(200,-4),(250,41)}
\put(180,60){\circle*{5}}
\put(220,60){\circle*{5}}

\put(120,91){\line(2,-1){60}}

\end{picture}
\caption{A Seifert graph with a cut vertex, whose knot diagram is flat.}
\label{FlatCutVertex}
\end{figure}

\begin{Theorem}
\label{block}
Diagrams of prime knots are flat if and only if their Seifert graphs are blocks.
\end{Theorem}
\textbf{Proof}
If the Seifert graph is a block then it has no cut vertices, and so from Lemma \ref{notblock} the knot diagram is flat. Conversely, consider a flat knot diagram, and suppose that the Seifert graph has a cut vertex. Now, by reconstructing that part of the knot diagram in the neighbourhood of the cut vertex, as shown in Figure \ref{cutvertex}, we find that the knot cannot be prime.

\setlength{\unitlength}{1pt}
\begin{figure}[h]
\begin{picture}(100,50)(-100,30)

\put(-80,30){\begin{picture}(100,100)(-70,40)

\put(10,60){\line(-6,1){60}}
\put(10,60){\line(-6,-1){60}}

\put(-50,65){\circle*{1}}
\put(-51,60){\circle*{1}}
\put(-50,55){\circle*{1}}

\put(10,60){\circle*{5}}

\put(10,60){\line(6,1){60}}
\put(10,60){\line(6,-1){60}}

\put(70,65){\circle*{1}}
\put(71,60){\circle*{1}}
\put(70,55){\circle*{1}}

\end{picture}}

\put(80,30){\begin{picture}(100,100)(-70,40)

\qbezier{(-30,65),(10,75),(50,65)}
\qbezier{(-30,55),(10,45),(50,55)}
\put(60,60){\circle{40}}
\put(57,57){$K_{2}$}
\put(-40,60){\circle{40}}
\put(-47,57){$K_{1}$}

\end{picture}}

\end{picture}

\caption{The neighbourhood of a cut vertex, and its knot diagram.}
\label{cutvertex}
\end{figure}
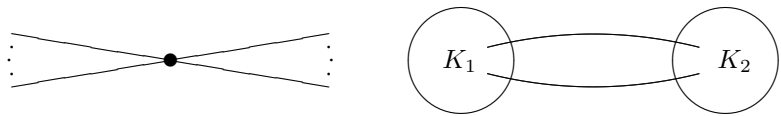

\begin{figure}[h]
\centering
\includegraphics[height=3cm]{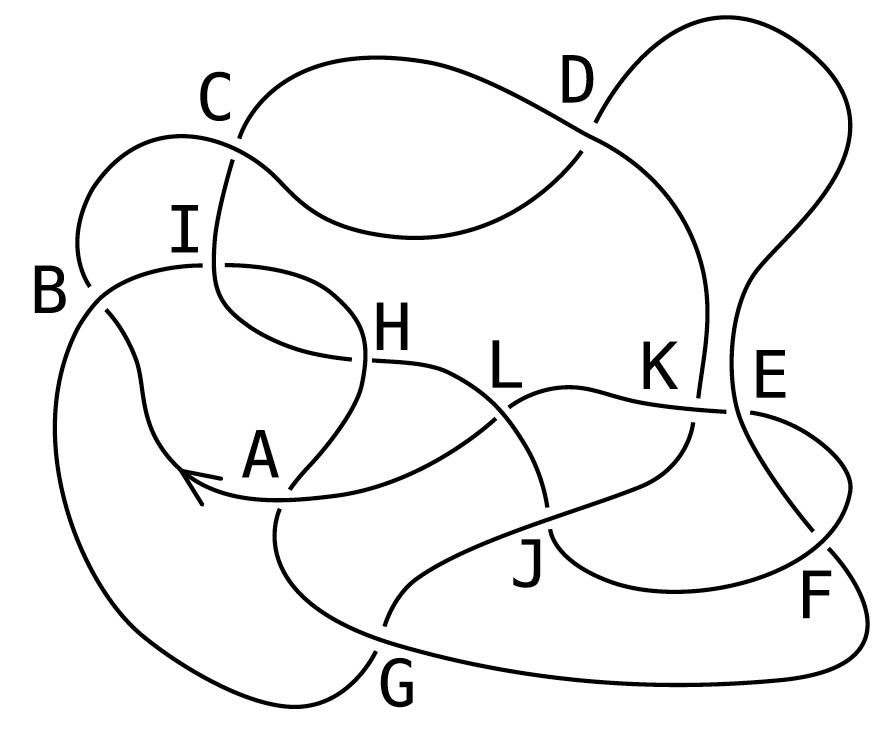}
\caption{The next smallest flat diagram after the trefoil.}
\label{12crossing}
\end{figure}

The usual diagram of the trefoil is flat. The next smallest flat diagram is shown in Figure \ref{12crossing}.

\begin{Theorem}
\label{flatimpliesST}
Flat knot diagrams have isomorphic Seifert and Tait graphs.
\end{Theorem}
\textbf{Proof}
Given a knot diagram $D$ we suppose that its Seifert graph $\Phi^{*}$ is not isomorphic to its Tait graph $T$. In this case, $T$ must have both c-edges and d-edges. Therefore, the graph $\Phi$ has more than one component, and this implies that $\Phi^{*}$ has a cut vertex. But then, from Theorem \ref{block}, $D$ is not flat. \hfill$\Box$

The converse is false. The knot $5_{1}$ is the first of many counterexamples.

\section{The cyclic word algorithm}

We start by describing how a link diagram leads to a graph cellularly embedded on an orientable surface.

Given a link diagram with $\ell$ components, label each crossing uniquely from the alphabet $a_{1},a_{2},\dots$. Then, for each component of the link, draw a polygon (topologically a disc) whose sides correspond to the crossings of that component, and transfer the labelling from the crossings to the sides. This yields a collection of $\ell$ polygons, one for each component. In this collection, each label appears twice, once for an over-crossing and again for the corresponding under-crossing. Now, for each over-crossing labelled $a_{i}$, say, re-label the under-crossing by $a_{i}^{-1}$. Finally, identify pairs of sides in the collection of polygons whenever they are labelled by a letter and its inverse. We obtain an orientable surface, and the identified sides become edges of a graph cellularly embedded in the surface. The labels around each face of this embedded graph form a cyclic word, so we have $\ell$ such words.

In the case of knots, $\ell=1$, our graph has just one face, which defines a single cyclic word. Note that in this case the labelled polygon is equivalent to the \emph{Gauss diagram} of the knot, which is constructed from the knot diagram by connecting by an oriented chord the pre-images of each double point of the immersion. The orientation of the chords goes from the over-passing branch to the under-passing one.

\begin{figure}[h]
\centering
\includegraphics[height=3cm]{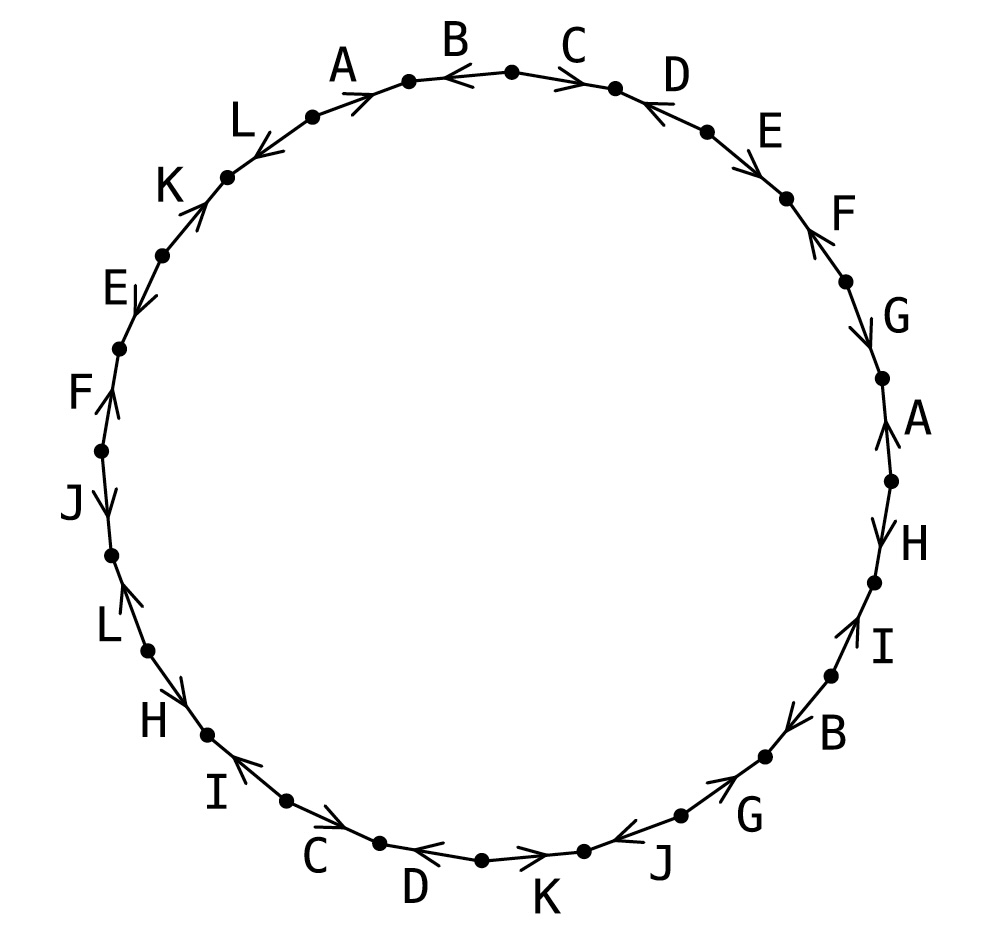}
\caption{The word corresponding to the diagram in Figure \ref{12crossing}.}
\label{12crossingword}
\end{figure}

\begin{Definition}[\cite{Alina2005}]
We call this algorithm for constructing a surface from a knot diagram the \emph{cyclic word algorithm}.
\end{Definition}
See Figure \ref{12crossingword}, for example.

Now restrict to \emph{prime} knots, and only consider their \emph{reduced} diagrams, that is, diagrams without any nugatory crossings. Then, the word defined as above is a \emph{Wicks form}:

\begin{Definition}[\cite{Bacher}]
\label{Wicks}
A \emph{Wicks form} is a cyclic word $w$ in some alphabet $a_{1},a_{2},\dots$ with inverses $a_{1}^{-1},a_{2}^{-1},\dots$ such that for $\epsilon, \delta\in\{-1,1\}$
\begin{description}
\item[\rm (i)] if $a_{i}^{\epsilon}$ appears in $w$ then $a_{i}^{-\epsilon}$ appears exactly once in $w$,
\item[\rm (ii)] $w$ contains no cyclic factor of the form $a_{i}^{\epsilon}a_{i}^{-\epsilon}$, and
\item[\rm (iii)] if $a_{i}^{\epsilon}a_{j}^{\delta}$ is a cyclic factor of $w$ then $a_{j}^{-\delta}a_{i}^{-\epsilon}$ is not a cyclic factor of $w$.
\end{description}
\end{Definition}

\begin{figure}[h]
\centering
\includegraphics[height=3cm]{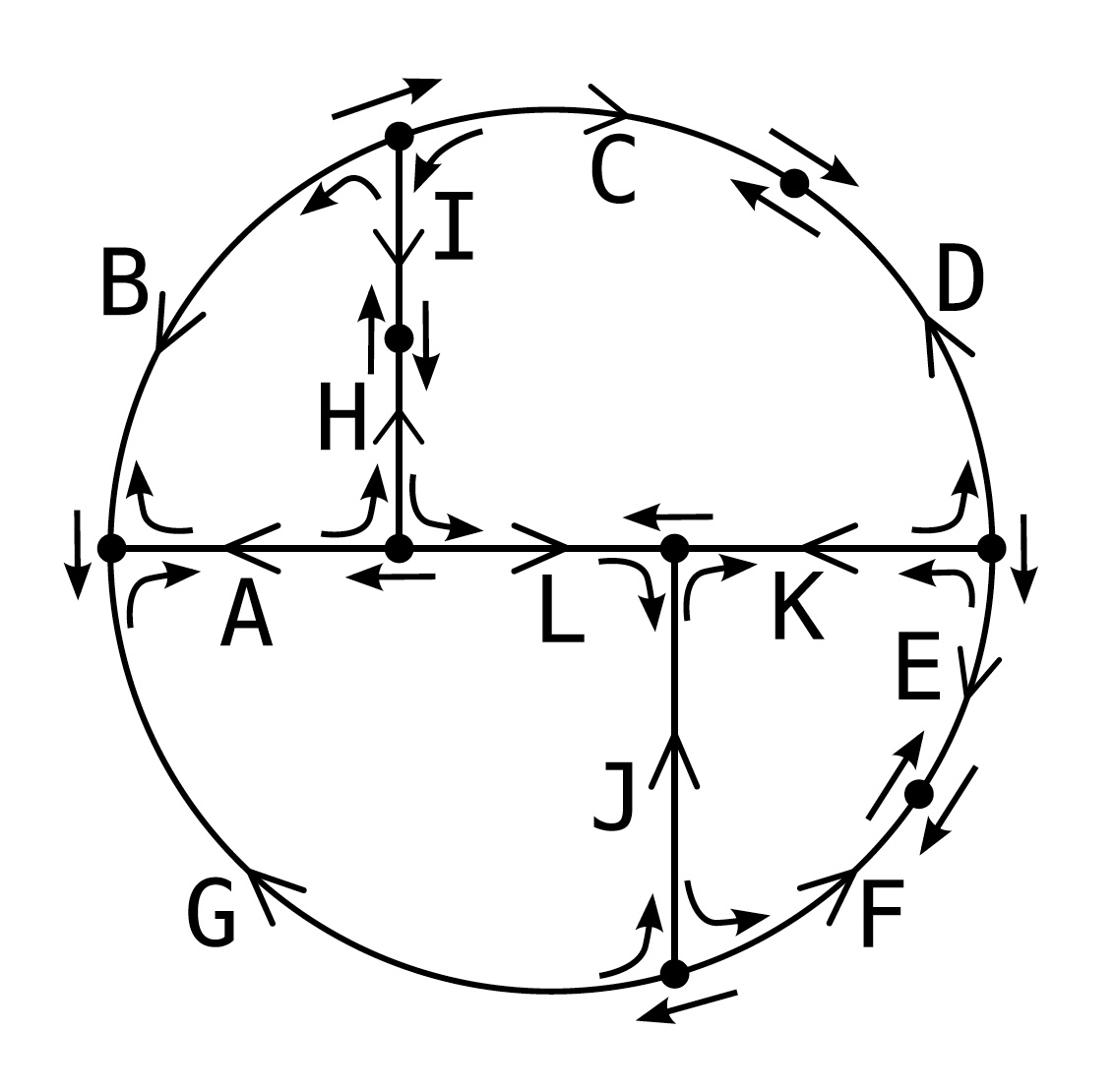}
\caption{The graph corresponding to the diagram in Figure \ref{12crossing}. It is bipartite, with no cut-vertices. All of the edges of the Tait graph are d-edges.}
\label{12crossinggraph}
\end{figure}

Suppose we start with a Wicks form $w$. Write $w$ around a $2e$-gon, one letter on each side. Then, following the algorithm above, for $i\in\{1,\dots,e\}$ identify the sides labelled $a_{i}$ and $a_{i}^{-1}$ of the $2e$-gon, obtaining the surface $S$. Also, record where the (identified) sides of the $2e$-gon are, to obtain a graph $\G$ (with $e$ edges and $v$ vertices) lying on $S$, having just one face. The Wicks form $w$ then determines a \emph{bieulerian path} on $G$:

\begin{Definition}[\cite{Alina2005}]
A \emph{bieulerian path} in a graph is a closed path that traverses each edge exactly twice, once in each direction, and does not traverse any edge followed immediately by its inverse.
\end{Definition}
For example, the word in Figure \ref{12crossingword} determines the bieulerian path in Figure \ref{12crossinggraph}.

\begin{Definition}
A graph with a bieulerian path is \emph{realizable} if it arises from a knot diagram via a Wicks form.
\end{Definition}

\begin{Definition}
Given a component of the boundary of a surface, we may identify it with the boundary of a disc. This is called \emph{capping} that boundary component.
\end{Definition}

\begin{Theorem}
Given a knot diagram, the surface obtained using Seifert’s algorithm and then capping it is homeomorphic to that obtained from the cyclic word algorithm.
\end{Theorem}

\textbf{Proof}
Given a knot $K$, we draw its Gauss diagram $G$. We will construct an orientable surface $F$ corresponding to $G$ (see page 218 of \cite{Alina2002}). We first ``thicken'" $G$, by replacing its circle by an annulus and each chord by a ribbon whose ends are glued to the ``inner" edge of the annulus in such a way that the resulting surface, which has a number of boundary components, is orientable. We cap each of the boundary components, except for one. The boundary component left uncapped is the ``outer" edge of the thickened circle of the original Gauss diagram $G$. We now have a surface we denote by $F'$, with one boundary component. $F'$ is the canonical Seifert surface of $K$.

Now cap the remaining boundary component of $F'$ with a disc, to obtain the required closed orientable surface $F$. On this surface the original Gauss diagram $G$ may be regarded as an embedded trivalent graph with $3e$ edges, $2e$ vertices, and $v+1$ faces (where $e$ and $v$ are as above in the description of the graph $\G$). Of these $3e$ edges, $2e$ are ``circle edges" and $e$ are ``chord" edges. We have
$$\chi(F)=v+1-3e+2e=v+1-e.$$

Next, return to the thickened Gauss diagram, and cap the outer edge of the annulus with a disc, leaving all the other boundary components uncapped. This gives a surface we denote by $S'$. Then, make the ribbons wider and wider until none of the inner edge of the annulus is left. Finally, shorten all the ribbons until they have length zero. This recovers the closed surface $S$ defined above from the cyclic word algorithm.

$F$ and $S$ are both orientable and $\chi(F)=\chi(S)$, so they are homeomorphic.

\hfill$\Box$

\begin{figure}[h]
\centering
\includegraphics[width=11cm, trim= 2cm 14cm 2cm 10cm, clip=true]{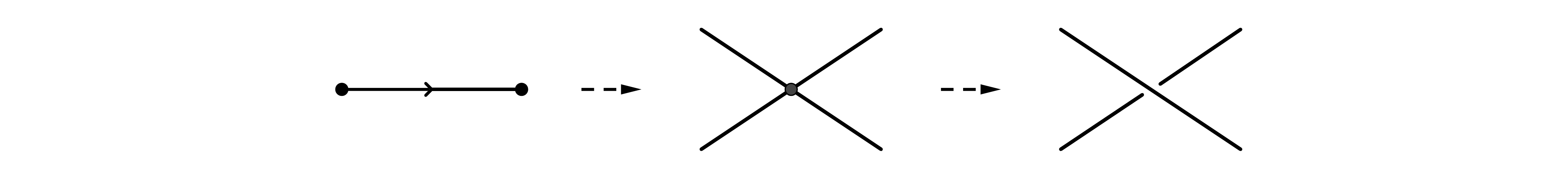}
\caption{From an oriented edge to a crossing.}
\label{Medial}
\end{figure}

Let $P$ be a planar graph. There is an algorithm, which is the inverse of Seifert’s algorithm, for constructing a link diagram from $P$, as follows. Firstly, we choose an arbitrary orientation for each of the vertices of the graph. Secondly, whenever an edge joins vertices with the same orientation we add a vertex of degree two, to obtain the graph $Q$. Then all edges of $Q$ will join vertices of opposite orientations, and we can orient these edges so that they point from an anti-clockwise to a clockwise vertex. Let $M$ be the medial graph of $Q$. Replace each vertex of $M$ with a crossing, determined as in Figure \ref{Medial}, to obtain a link diagram.

\begin{Theorem}[4.3 of \cite{Alina2005}]
\label{planar}
A planar trivalent graph with a bieulerian path corresponds to a flat knot diagram.
\end{Theorem}
\textbf{Proof}
Suppose that $G$ is a planar trivalent graph with a bieulerian path, and consider a vertex $p$ in $G$. When the path first arrives at $p$ it must leave it on one of the other two edges, either to the left or to the right. This choice induces an orientation on the vertex $p$, and hence on all the vertices, and we can now use the algorithm to construct a link diagram. But walking around the bieulerian path in $G$ is the same as following a component of $D$ and meeting every crossing twice. So $D$ has just one component: it is a knot diagram. By construction the Seifert circles are all of type I, and so it is a flat diagram. \hfill$\Box$

\section{Seifert-Tait knots dominate}

\begin{figure}[h]
\centering
\includegraphics[height=1.5cm]{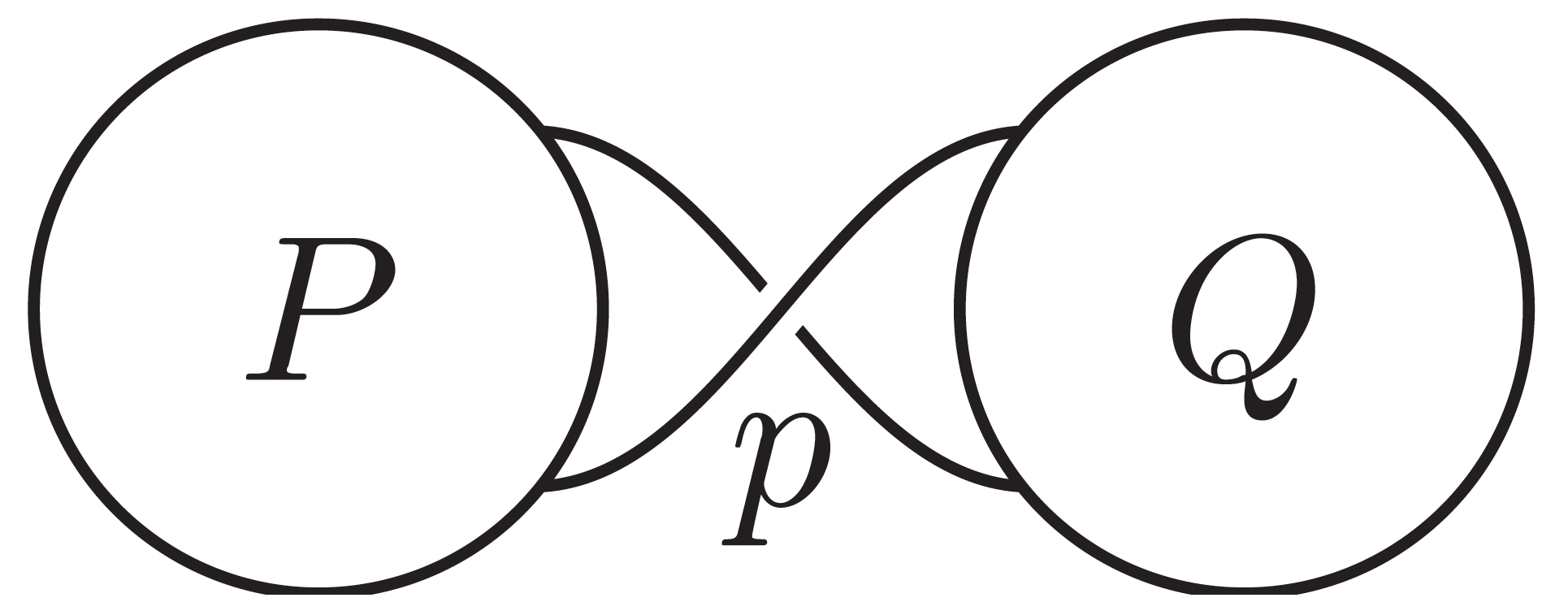}
\caption{$p$ is a nugatory crossing.}
\label{Nugatory}
\end{figure}

We start by recalling some elementary properties of Gauss diagrams \cite{Alina2002}. Let $D$ be an alternating diagram of a knot $K$, and let $G$ be its Gauss diagram and $w$ the corresponding cyclic word. Suppose a chord in $G$ labelled by the letters $a$ and $a^{-1}$ in $w$ is such that there is no other chord $(b,b^{-1})$ appearing in $w$ as
\begin{equation}
\dots a\dots b^{-1}\dots a^{-1}\dots b
\end{equation}
(which might be because $a$ and $a^{-1}$ are adjacent). Then the corresponding crossing $p$ in $D$ is \emph{nugatory}, and the chord in $G$ is \emph{isolated}: see Figure \ref{Nugatory}.

If $K$ is prime and $D$ is a reduced diagram, then $G$ cannot have any isolated chords.

\begin{Definition}
A pair of chords $(a,a^{-1})$ and $(b,b^{-1})$ in $G$ is called \emph{parallel} if their endpoints occur as
$$\dots ab^{-1}\dots ba^{-1}\dots$$
in $w$.
\end{Definition}
This corresponds to a \emph{clasp} in $D$: see Figure \ref{Clasp}.

\begin{figure}[h]

\setlength{\unitlength}{.7mm}
\thicklines

\begin{picture}(60,20)(40,10)

\qbezier{(105,25),(110,25),(114,21)}
\qbezier{(125,15),(120,15),(116,19)}

\qbezier{(105,15),(110,15),(115,20)}
\qbezier{(115,20),(120,25),(125,25)}

\qbezier{(125,25),(130,25),(134,21)}
\qbezier{(145,15),(140,15),(136,19)}

\qbezier{(125,15),(130,15),(135,20)}
\qbezier{(135,20),(140,25),(145,25)}
    
\end{picture}
\caption{A clasp}
\label{Clasp}
\end{figure}

\begin{Definition}
A triple of chords in $G$ is called \emph{parallel} if both the first two and the last two are parallel pairs.
\end{Definition}

\begin{figure}[h]
\centering
\includegraphics[height=1.5cm]{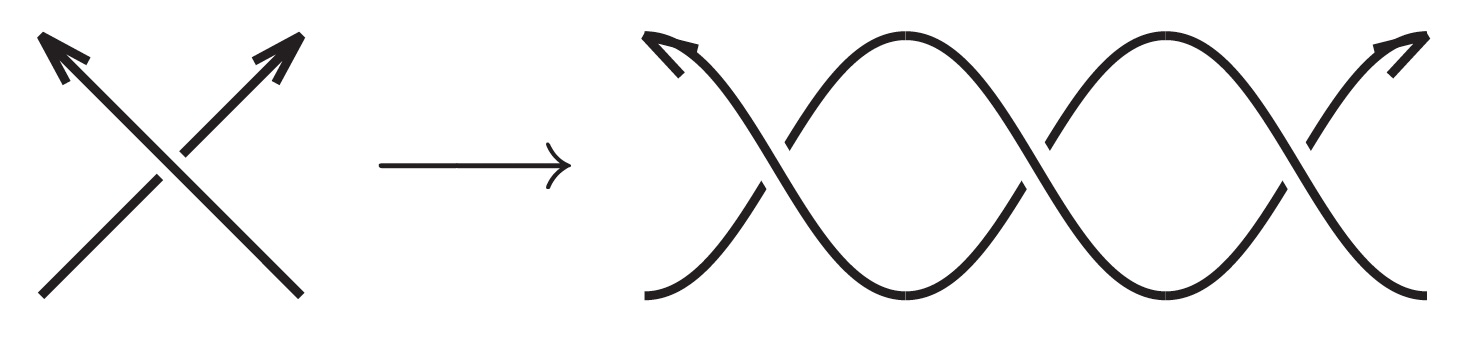}
\caption{The $\bar t_{2}$ move.}
\label{t2move}
\end{figure}

\begin{Definition}
The move illustrated in Figure \ref{t2move} is called the \emph{$\bar t_{2}$ move}.
\end{Definition}

\begin{Definition}
A knot diagram is called \emph{$\bar t_{2}$-irreducible} if its Gauss diagram does not have a triple of parallel chords.
\end{Definition}

The importance of this definition is that for a given genus $g$ we can divide all the alternating diagrams into a finite number of sets, called \emph{generating series}, as follows.

\begin{Theorem}[\cite{Stoimenow}]
The set of $\bar t_{2}$-irreducible alternating diagrams with no nugatory crossings and genus $g$ is finite.
\end{Theorem}
\begin{Definition}
Elements of this set are called \emph{generating diagrams}, and the set of all diagrams obtained from a single generating diagram using $\bar t_{2}$ moves is called a \emph{generating series}.
\end{Definition}

Next, we introduce a key result on alternating diagrams using \emph{flypes} which was proved by Menasco and Thistlethwaite in \cite{Menasco}.

\begin{figure}[h]
\centering
\includegraphics[height=2cm]{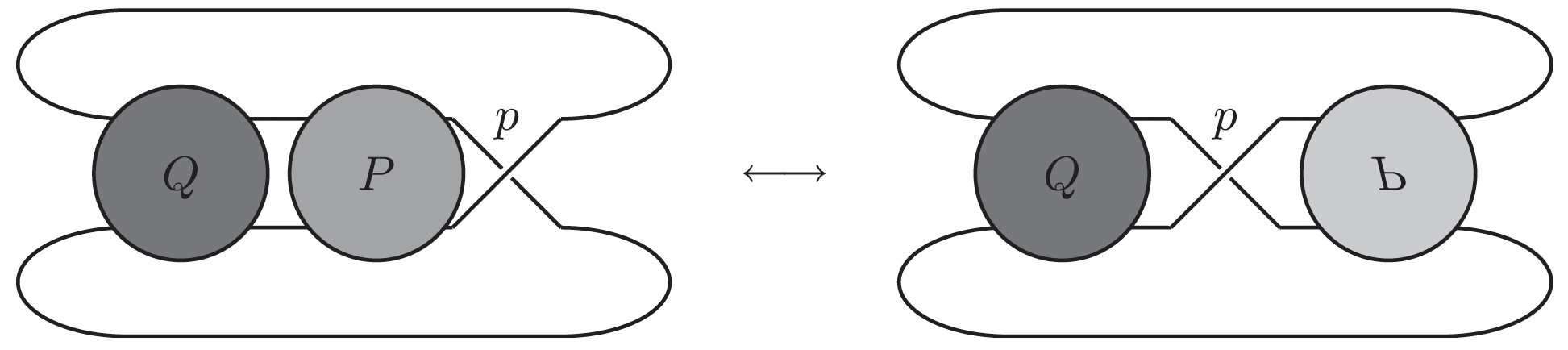}
\caption{A flype near the crossing $p$.}
\label{Flype}
\end{figure}

\begin{Definition}
A \emph{flype} is a move on a diagram illustrated in Figure \ref{Flype}.
\end{Definition}

\begin{Theorem}[\cite{Menasco}]
Any two reduced, prime, oriented, alternating diagrams of the same link are related by a sequence of flypes.
\end{Theorem}
\begin{Corollary}
\label{unique}
If an alternating diagram admits no flypes then it is the only alternating diagram of its knot.
\end{Corollary}

In order to be able to make use of Corollary (\ref{unique}) we need the following equivalence relation on the set of crossings in a diagram.

\begin{Definition}
Two crossings $a$ and $b$ in a diagram $D$ are \emph{$\sim$-equivalent} if their chords in the Gauss diagram can be made parallel by a sequence of flypes on $D$.
\end{Definition}

\begin{Theorem}[2.14 and 4.2 of \cite{Alina2005}]
A knot diagram of genus $g$ has at most $6g-3$ \mbox{$\sim$-equivalence} classes, and this bound is attained.
\end{Theorem}

\begin{Theorem}[4.7 of \cite{Alina2005}]
\label{noflypes}
Let $G$ be a planar trivalent graph with a bieulerian path, and $D$ its knot diagram. Then the following conditions are equivalent:
\begin{enumerate}
\item $G$ is three-connected,
\item $D$ has $6g-3$ $\sim$-equivalence classes,
\item $D$ admits no (non-trivial) flypes.
\end{enumerate}
\end{Theorem}
The diagram $D$ in this result must be flat, by Theorem \ref{planar}. Also, by Corollary (\ref{unique}), it must be the only alternating diagram of its knot. Therefore, three-connected planar trivalent graphs with a bieulerian path are complete invariants for flat alternating knots, and hence for Seifert-Tait knots.

\begin{Theorem}[\cite{Alina2005}]
\label{limit}
For a fixed genus $g>1$,
$$\frac{\#\{\hbox{flat alternating knots with crossing number }n\hbox{ and genus }g\}}{\#\{\hbox{prime alternating knots with crossing number }n\hbox{ and genus }g\}}\longrightarrow 1$$
as $n\longrightarrow\infty$.
\end{Theorem}
\textbf{Proof}
For alternating knots of genus $g$, the maximum number of $\sim$-equivalence classes is $6g-3$, and this maximum is attained. Consider a generating series in which the knots have $6g-3$ $\sim$-equivalence classes, and are therefore flat. Knots in this series with $n$ crossings arise from partitions of $n$ into the $6g-3$ classes, and so there are
$$\binom{n+6g-3-1}{6g-3-1}=\frac{(n+6g-4)!}{(6g-4)!n!}$$
of them.

Now consider the generating series which have fewer than $6g-3$ $\sim$-equivalence classes. Suppose there are $c$ of them, and consider one with $d$ $\sim$-equivalence classes. There will be
$$\binom{n+d-1}{d-1}$$
knots in this series.

But $d\leq 6g-4$, so the number of knots in all of the non-flat series taken together must be no more than
$$c\binom{n+6g-4-1}{6g-4-1}=\frac{c(n+6g-5)!}{(6g-5)!n!}.$$

Therefore the required ratio is at least
$$\frac{(n+6g-4)!}{(6g-4)!n!}\times\frac{(6g-5)!n!}{c(n+6g-5)!}=\frac{(n+6g-4)}{c(6g-4)}$$
and hence the result.\hfill$\Box$

\begin{Theorem}
\label{STresult}
For a fixed genus $g>1$,
$$\frac{\#\{\hbox{Seifert-Tait alternating knots with crossing number }n\hbox{ and genus }g\}}{\#\{\hbox{prime alternating knots with crossing number }n\hbox{ and genus }g\}}\longrightarrow 1$$
as $n\longrightarrow\infty$.
\end{Theorem}
\textbf{Proof}
This follows from Theorems \ref{flatimpliesST} and \ref{limit}.\hfill$\Box$

\bigskip

\noindent
\textsc{Stephen Huggett}\\
Centre for Mathematical Sciences, University of Plymouth, UK.

\bigskip

\noindent
\textsc{Alina Vdovina}\\
City University of New York, City College and Graduate Center, New York,~USA.


\begin{thebibliography}{99}

\bibitem{Bacher}R.~Bacher, A.~Vdovina \emph{Counting 1-vertex triangulations of oriented surfaces} Discrete Mathematics \textbf{246}(1--3) 13--27, 2002.

\bibitem{Cromwell1989}P.~Cromwell \emph{Homogeneous Links} J. London Math. Soc. (2) \textbf{39}, 535--552, 1989.

\bibitem{Cromwell2004}P.~Cromwell \emph{Knots and Links} Cambridge Univ. Press, Cambridge, 2004.

\bibitem{Dasbachetal}O.~Dasbach, D.~Futer, E.~Kalfagianni, X-S.~Lin, and N.~Stoltzfus \emph{The Jones polynomial and graphs on surfaces} Journal of Combinatorial Theory, Series B \textbf{98}, 384--399, 2008.

\bibitem{HMV}S.~Huggett, I.~Moffatt, N.~Virdee \emph{On the Seifert graphs of a link diagram and its parallels} Math. Proc. Camb. Phil. Soc. \textbf{153}, 123--145, 2012.

\bibitem{Bipartite}S.~Huggett and I.~Moffatt \emph{Bipartite partial duals and circuits in medial graphs} Combinatorica, \textbf{33}(2), 231--252, 2013.

\bibitem{Kharlampovich}O.~Kharlampovich, A.~Vdovina \emph{Low complexity algorithms in knot theory} Int. J. Alg. Comp. \textbf{29}(2), 245--262, 2019.

\bibitem{Menasco}W.~Menasco, M.~Thistlethwaite \emph{The classification of alternating links} Annals of Mathematics \textbf{138}, 113--171, 1993.

\bibitem{Murasugi}K.~Murasugi \emph{On alternating knots} Osaka Math. J. \textbf{12}, 277--303, 1960.

\bibitem{Alina2005}A.~Stoimenow, A.~Vdovina \emph{Counting alternating knots by genus} Math. Ann. \textbf{333}, 1--27, 2005.

\bibitem{Alina2002}A.~Stoimenow, V.~Tchernov, A.~Vdovina \emph{The canonical genus of a classical and virtual knot}
Geom. Dedicata \textbf{95}, 215--225, 2002.

\bibitem{Stoimenow}A.~Stoimenow \emph{Knots of genus one} Proc. Amer. Math. Soc. \textbf{129}(7), 2141--2156, 2001.

\bibitem{Turaev}V.~Turaev \emph{A simple proof of the Murasugi and Kauffman theorems on alternating links} Enseign. Math. (2) \textbf{33}(3--4) 203--225, 1987.

\end{thebibliography}
\end{document}